\documentclass{article}%
\usepackage{amssymb}
\usepackage{graphicx}
\usepackage{amsmath}
\usepackage{amsfonts}%
\providecommand{\U}[1]{\protect\rule{.1in}{.1in}}
\newtheorem{theorem}{Theorem}

\newtheorem{corollary}[theorem]{Corollary}

\newtheorem{lemma}[theorem]{Lemma}

\newtheorem{proposition}[theorem]{Proposition}
\newtheorem{remark}[theorem]{Remark}

\begin{document}

\begin{center}
{\large On the local limit theorems for lower psi-mixing Markov chains.}

\end{center}

\bigskip Florence Merlev\`{e}de,\ Magda Peligrad and Costel Peligrad

\begin{center}
\bigskip
\end{center}

Universit\'{e} Gustave Eiffel

LAMA and CNRS UMR 8050.

5, Boulevard Descartes, Champs sur Marne

Marne La Vall\'{e}e, C\'{e}dex 2 France

E-mail address: florence.merlevede@univ-eiffel.fr

\bigskip

Department of Mathematical Sciences

University of Cincinnati,

POBox 210025, Cincinnati, Oh 45221-0025, USA. \texttt{ }

E-mail address: peligrm@ucmail.uc.edu, peligrc@ucmail.uc.edu

\bigskip

\noindent\textbf{Keywords}: \noindent Markov chains, local limit theorem,
mixing, stable laws.

\bigskip


\begin{center}
\bigskip\bigskip

Abstract
\end{center}

In this paper we investigate the local limit theorem for additive functionals
of nonstationary Markov chains that converge in distribution. We consider both
the lattice and the non-lattice cases. The results are also new in the
stationary setting and lead to local limit theorems linked to convergence to
stable distributions. The conditions are imposed to individual summands and
are expressed in terms of lower psi-mixing coefficients.

\section{Introduction}

The local limit theorem for partial sums $(S_{n})_{n\geq1}$ of a sequence of
random variables deals with the rate of convergence of the probabilities of
the type $\mathbb{P}(c\leq S_{n}\leq d)$. Local limit theorems have been
studied for\ the case of lattice random variables and the case of non-lattice
random variables. A random variable is said to have a lattice distribution if
there exists $h>0$ and $\ell\in\mathbb{R}$ such that its values are
concentrated on the lattice $\{\ell+kh:k\in\mathbb{Z}\}$. The non-lattice
distribution means that no such $\ell$ and $h$ exist.

This problem was intensively studied for sums of i.i.d. random variables. For
i.i.d. random variables in the domain of attraction of a stable law the local
limit theorem was solved by Stone (1965) and Feller (1967).

It should be mentioned that the local limit theorem is more delicate than its
convergence in distribution counterpart and often requires additional
conditions. An important counterexample is given by Gamkrelidze (1964),
pointing out this phenomenon for independent summands, and a variety of
sufficient conditions were developed over the years. We mention especially
papers by Rozanov (1957), Mineka-Silverman (1970), Maller (1978), Shore (1978)
and Dolgopyat (2016).

In the dependent case, we mention early works on Markov chains by Kolmogorov
(1962). On the other hand, Nagaev (1961) gives rates of convergence in the CLT
for stationary $\psi-$mixing Markov chains. In the lattice case, Szewczak
(2010) established a local limit theorem for continued fractions, which is an
example of $\psi-$mixing sequence. Hafouta and Kifer (2016) proved a local
limit theorem for nonconventional sums of stationary $\psi-$mixing Markov
chains, while the case of infinite variance is analyzed in Aaronson and Denker
(2001a). Also in the stationary case we mention the local limit theorems for
Markov chains in the papers by Herv\'{e} and P\`{e}ne (2010), Ferr\'{e} et al.
(2012). Recently, there are two additional works by Merlev\`{e}de et al.
(2021) and Dolgopyat and Sarig (2021), which cover several aspects of local
limit theorems associated with the central limit theorem\ for $\psi-$mixing
nonstationary Markov chains.

These results raise the natural question if a local limit theorem is valid for
more general Markov chains. In this paper we positively answer this question
and consider a class larger than $\psi-$mixing Markov chains defined by using
the so called one sided lower $\psi-$mixing coefficient. The key tool in
proving these results is a delicate factorization of the characteristic
function of partial sums.

We shall comment that all the results obtained by Merlev\`{e}de et al. (2021)
for $\psi-$mixing nonstationary Markov chains are also valid for this larger
class. Furthermore, we shall also obtain local limit theorems associated with
other limiting distributions than the normal attraction, which is actually the
attraction to a stable distribution with index $p=2$. More precisely we shall
also obtain the local limit theorem associated to attraction to any stable
distributions with index $1<p<2$. The key tools for obtaining these results
are some general local limit theorems, which assume convergence to
distributions with integrable characteristic functions.

Our paper is organized as follows. In Section 2 we present some preliminary
considerations about the types of local limit theorems we shall obtain. In
Section 3 we present the main results and provide some examples of classes of
lower $\psi-$mixing Markov chains. Section 5 contains the proofs. We present
first two general local limit theorems, then a bound on the characteristic
function of sums and the proofs of the main results.

In the paper, by $\Rightarrow$ we denote the convergence in distribution.

\section{Preliminary considerations}

We formulate first the conclusions of the local limit theorems. Let
$(S_{n})_{n}$ be a sequence of random variables and let $L$\ be a random
variable with characteristic function $f_{L}$. The basic assumption will be an
underlying convergence in distribution:%
\begin{equation}
\frac{S_{n}}{B_{n}}\Rightarrow L,\text{ where }f_{L}\text{ is integrable and
}B_{n}\rightarrow\infty. \label{conv distribution}%
\end{equation}

Note that assuming that $f_{L}$ is integrable implies that $L$ has a
continuous density we shall denote by $h_{L}$ (see pages 370-371 in
Billingsley, 1999).

In the case when the variables $(S_{n})$ do not have values in a fixed minimal
lattice we shall say that the sequence satisfies a local limit theorem if for
any function $g$ on $\mathbb{R}$ which is continuous and with compact
support,
\begin{equation}
\lim_{n\rightarrow\infty}\sup_{u\in\mathbb{R}}\left\vert B_{n}{\mathbb{E}%
}g(S_{n}+u)-h_{L}\Big (\!\!-\frac{u}{B_{n}}\Big )\int g(t)\lambda
(dt)\right\vert =0,\label{GLLTnonlat}%
\end{equation}
where $\lambda$ is the Lebesgue measure.

If all the variables have values in the same fixed lattice
\[
\mathcal{S}=\{kh,k\in{\mathbb{Z}}\},
\]
where $h$ is some fixed positive number, by the local limit theorem we
understand that for any function $g$ on $\mathbb{R}$ which is continuous and
with compact support,
\begin{equation}
\lim_{n\rightarrow\infty}\sup_{u\in\mathcal{S}}\left\vert B_{n}{\mathbb{E}%
}g(S_{n}+u)-h_{L}\Big (\!\!-\frac{u}{B_{n}}\Big )\int g(v)\mathcal{L}%
_{h}(dv)\right\vert =0,\label{GLLTlatt}%
\end{equation}
where $\mathcal{L}_{h}$ is the measure assigning $h$ to each point $kh$.

By standard arguments, for any real numbers $c$ and $d$ such that $c\leq d$,
we also have that (\ref{GLLTnonlat}) implies in the non-lattice case that%
\[
\lim_{n\rightarrow\infty}\sup_{u\in\mathbb{R}}\left\vert B_{n}\mathbb{P}%
(c-u\leq S_{n}\leq d-u)-(d-c)h_{L}\Big (\!\!-\frac{u}{B_{n}}\Big )\right\vert
=0,
\]
and in the lattice case
\[
\lim_{n\rightarrow\infty}\sup_{u\in\mathcal{S}}\left\vert B_{n}\mathbb{P}%
(c-u\leq S_{n}\leq d-u)-\ h\sum\nolimits_{k}I(c-u\leq kh\leq d-u)h_{L}%
\Big (\!\!-\frac{u}{B_{n}}\Big )\right\vert =0.
\]
In particular, since $B_{n}\rightarrow\infty$ as $n\rightarrow\infty$, then
for fixed $A>0$, in the non-lattice case
\begin{equation}
\lim_{n\rightarrow\infty}\sup_{|u|\leq A}\left\vert B_{n}\mathbb{P}(c-u\leq
S_{n}\leq d-u)-(d-c)h_{L}(0)\right\vert =0, \label{Mallerresult}%
\end{equation}
and in the lattice case
\begin{equation}
\lim_{n\rightarrow\infty}\sup_{|u|\leq A,u\in\mathcal{S}}\left\vert
B_{n}\mathbb{P}(c-u\leq S_{n}\leq d-u)-\ h\sum\nolimits_{k}I(c-u\leq kh\leq
d-u)h_{L}(0)\right\vert =0. \label{Mlattice}%
\end{equation}
If we further take $u=0$ in (\ref{Mallerresult}), then, in the nonlattice case
we have
\[
\lim_{n\rightarrow\infty}B_{n}\mathbb{P}(S_{n}\in\lbrack c,d])=\left(
d-c\right)  h_{L}(0).
\]
In other words, the sequence of measures $B_{n}\mathbb{P}(S_{n}\in\lbrack
c,d])$ of the interval $[c,d]$ converges to a scalar multiple of the Lebesgue measure.

In the lattice case we can take $u=0$ in (\ref{Mlattice})\ and obtain,%
\[
\lim_{n\rightarrow\infty}B_{n}\mathbb{P}(S_{n}\in\lbrack c,d])=h\sum
\nolimits_{k}I(c\leq kh\leq d)h_{L}(0).
\]
\ \ \ We assume now that $(\xi_{k})_{k\geq0}$ is a Markov chain defined on
$(\Omega,\mathcal{K},\mathbb{P})$ with values in $(\mathcal{X},\mathcal{B}%
(\mathcal{X})),$ where $\mathcal{B}(\mathcal{X})$ is a $\sigma-$field on
$\mathcal{X}$, with regular transition probabilities:
\begin{equation}
Q_{k}(x,A)=\mathbb{P}(\xi_{k}\in A|\xi_{k-1}=x)\label{defQ}%
\end{equation}
and marginal distributions denoted by%
\begin{equation}
\mathbb{P}_{k}(A)=\mathbb{P}(\xi_{k}\in A).\label{def Marginal}%
\end{equation}
Throughout the paper we shall assume that there is a constant $a>0$ with the
following property:

For all $k\geq1$ there is $\mathcal{X}_{k}^{\prime}\in\mathcal{B}%
(\mathcal{X})$ with $\mathbb{P}_{k-1}(\mathcal{X}_{k}^{\prime})=1$ such that
for all $A\mathcal{\in B}(\mathcal{X})$ and $x\in\mathcal{X}_{k}^{\prime}$ we
have%
\begin{equation}
Q_{k}(x,A)\geq a\mathbb{P}_{k}(A). \label{main condition}%
\end{equation}
Triangular arrays $(\xi_{kn})_{k\geq0}$ can also be considered. In this case
$a$ should be replaced by $a_{n}$ which may be different from line to line.

Let $(g_{j})_{j\geq0}$ be real-valued measurable functions and define
\begin{equation}
X_{j}=g_{j}(\xi_{j}) \label{defX}%
\end{equation}
and set
\[
S_{n}=\sum_{k=1}^{n}X_{k}.
\]

In order to obtain our results we shall combine several techniques,
specifically designed for obtaining local limit theorems, with a suitable
bound on the characteristic function facilitated by condition
(\ref{main condition}).

\section{Main Results}

In the sequel, we shall denote by $f_{k}(t)$ the characteristic function of
$X_{k},$%
\[
f_{k}(t)=\mathbb{E}(\exp(itX_{k})).
\]

We shall introduce two useful conditions that will provide a version of the
local limit theorem under conditions imposed to the characteristic function.

\smallskip

Let $a$ be as in (\ref{main condition}).

\smallskip

\textbf{Condition A}. There is $\delta>0$ and $n_{0}\in{\mathbb{N}}$ and a
Borel function $g:[1,\infty)\rightarrow(0,\infty)$ such that for $1\leq
|u|\leq\delta B_{n}$ and $n>n_{0}$
\begin{equation}
\frac{a^{2}}{2^{4}}\sum_{k=1}^{n} \Big (1- \big |f_{k} \big (\frac{u}{B_{n}%
}\big ) \big |^{2} \Big )>g(|u|)\text{ and }\exp(-g(|u|))\text{ is
integrable.} \label{regular}%
\end{equation}

\textbf{Condition B}. For $u\neq0$ there is an $\varepsilon=\varepsilon(u)$,
$c(u)$ and a $n_{0}=n_{0}(u)$ such that for all $t$ with $|t-u|\leq
\varepsilon$ and $n>n_{0}$
\begin{equation}
\frac{a^{4}}{2^{4}(\ln B_{n})}\sum_{k=1}^{n}(1-|f_{k}(t)|^{2})\geq c(u)>1.
\label{verC2}%
\end{equation}
Our general local limit theorem is the following:

\begin{theorem}
\label{ThLocalAB}Let $(X_{j})_{j\geq0}$ be defined by (\ref{defX}) and
satisfying (\ref{main condition}). Assume that conditions
(\ref{conv distribution}), \textbf{A} and \textbf{B} are satisfied. If not all
the variables have values in a fixed lattice, then (\ref{GLLTnonlat}) holds.
If all the variables have the values in a fixed lattice $\mathcal{S}%
=\{kh,k\in{\mathbb{Z}}\}$, under the same conditions with the exception that
we assume now that Condition\textbf{ B} holds for $0<|u|\leq\pi/h,$ then
(\ref{GLLTlatt}) holds.
\end{theorem}

The conditions of Theorem \ref{ThLocalAB} are verified in many situations of
interest. Merlev\`{e}de et al. (2021) provided sufficient conditions for
Conditions \textbf{A} and \textbf{B} when second moment is finite and obtained
the local limit theorems in the nonlattice case under a more restrictive
condition than (\ref{main condition}). Namely it was assumed there that there
are two constants $a>0$ and $b<\infty$ with the following property:

For all $k\geq1$ there is $\mathcal{X}_{k}^{\prime}\in\mathcal{B}%
(\mathcal{X})$ with $\mathbb{P}_{k-1}(\mathcal{X}_{k}^{\prime})=1$ such that
for all $A\mathcal{\in B}(\mathcal{X})$ and $x\in\mathcal{X}_{k}^{\prime}$ we
have
\begin{equation}
a\mathbb{P}_{k}(A)\leq Q_{k}(x,A)\leq b\mathbb{P}_{k}(A)\text{.} \label{cond0}%
\end{equation}
By using our Proposition \ref{factor}, all the results in Merlev\`{e}de et al.
(2021) also hold for a larger class of Markov chains satisfying only the one
sided condition (\ref{main condition}). Therefore in all that results in
Merlev\`{e}de et al. (2021) the condition $b<\infty$ is superfluous.
Furthermore, because of Theorem \ref{ThLattice} below, all the results in
Merlev\`{e}de et al. (2021) can also be formulated for the lattice case. We
shall not repeat all the results there and mention only that our results are
also new in the stationary setting. In this particular case the results have a
simple formulation.

The first theorem deals with local limit theorem in case of attraction to
normal distribution, which is stable with index $p=2.$

\begin{theorem}
\label{normal attr}Assume that $(\xi_{k})_{k\in\mathbb{Z}}$ is a strictly
stationary Markov chain satisfying (\ref{main condition}). Define
$(X_{k})_{k\in\mathbb{Z}}$ by $X_{k}=g(\xi_{k})$ and assume ${\mathbb{E}%
}(X_{0})=0$ and $H(x)={\mathbb{E}}(X_{0}^{2}I(|X_{0}|\leq x))$ is a slowly
varying function as $x\rightarrow\infty$. Then, there is $B_{n}\rightarrow
\infty$ such that in the nonlattice case (\ref{GLLTnonlat}) holds and in the
lattice case (\ref{GLLTlatt}) holds with $h_{L}(x)=(2\pi)^{-1/2}\exp
(-x^{2}/2)$.
\end{theorem}

\begin{remark}
\label{remnormal attr}In the case where ${\mathbb{E}}(X_{0}^{2})<\infty,$ we
can take $B_{n}^{2}=\mathbb{E(}S_{n}^{2}),$ and there are two constants
$C_{1}$ and $C_{2}$ such that $C_{1}{\mathbb{E}}(X_{0}^{2})n\leq
\mathbb{E(}S_{n}^{2})\leq C_{2}{\mathbb{E}}(X_{0}^{2})n.$ Furthermore, since
in this case $\mathbb{E(}S_{n}^{2})/n\rightarrow c^{2}>0,$ we can also take
$B_{n}^{2}=c^{2}n$. In case ${\mathbb{E}}(X_{0}^{2})=\infty$ we can take
$B_{n}^{2}=\sqrt{\pi/2}\mathbb{E}|S_{n}|.$
\end{remark}

We move now to other types of limiting distributions $L$ in
(\ref{conv distribution}).

Let $0<p<2.$ Assume that $(\xi_{k})_{k\in\mathbb{Z}}$ is a strictly stationary
Markov chain and let $(X_{k})$ be defined by $X_{k}=g(\xi_{k}),$ with a
nondegenerate marginal distribution satisfying the following condition:
\begin{equation}
{\mathbb{P}}(|X_{0}|>x)=x^{-p}\ell(x), \label{tp}%
\end{equation}
where $\ell(x)$ is a slowly varying function at $\infty$,
\begin{equation}
\frac{{\mathbb{P}}(X_{0}>x)}{{\mathbb{P}}(|X_{0}|>x)}\rightarrow
c^{+}\;\;\;\text{and}\;\;\;\frac{{\mathbb{P}}(X_{0}<-x)}{{\mathbb{P}}%
(|X_{0}|>x)}\rightarrow c^{-}\text{ as }x\rightarrow\infty\label{lrtp}%
\end{equation}
with $0\leq c^{+}\leq1$ and $c^{+}+c^{-}=1.$

These conditions practically mean that the distribution of $X_{0}$ is in the
domain of attraction of a stable law with index $p.$ The convergence in
distribution of a mixing sequence of random variables in the domain of
attraction of a stable law is a delicate problem, which was often studied in
the literature, starting with papers by Davis (1983), Jakubowski (1991), Samur
(1987), Denker and Jakubovski (1989), Kobus (1995), Tyran-Kami\'{n}ska (2010),
Cattiaux and Manou-Abi (2014), El Machkouri et al. (2020), among many others.
In general an additional "non clustering" assumption is relevant for this type
of convergence. We shall assume that for every $x>0$ and for all $k\geq1$
\begin{equation}
\lim_{n\rightarrow\infty}\mathbb{P}(|X_{k}|>xB_{n}\text{ }|\text{ }%
|X_{0}|>xB_{n})=0, \label{DJ}%
\end{equation}
where $B_{n}$ is such that%
\begin{equation}
\frac{n}{B_{n}^{p}}\ell(B_{n})=n\mathbb{P}(|X_{0}|>B_{n})\rightarrow1.
\label{defBn}%
\end{equation}
In the context of $\varphi-$mixing sequences, for this type of condition we
refer to Theorem 3.2 in Samur (1997), Corollary 5.9 of Kobus (1995) and
Corollary 1.3 in Tyran-Kami\'{n}ska (2010), where the necessity of the
condition (\ref{DJ})\ for such a type of result is also discussed.

\begin{theorem}
\label{Th attraction}Assume $0<p<2$ and (\ref{tp}) and (\ref{lrtp}) are
satisfied. When $1<p<2$ we assume $\mathbb{E}(X_{0})=0$ and when $p=1$ we
assume $X_{0}$ has a symmetric distribution. Also assume (\ref{main condition}%
) and (\ref{DJ}). Then in the nonlattice case (\ref{GLLTnonlat}) holds and in
the lattice case (\ref{GLLTlatt}) holds with $B_{n}$ defined by (\ref{defBn}%
)$\ $and $h_{L}\ $the density of a strictly stable distribution $L$ with index
$p.$
\end{theorem}

Note that condition (\ref{DJ})\ of Theorem \ref{Th attraction} is satisfied if
we assume that there is $b<\infty$ such that the right hand side of
(\ref{cond0}) holds. Therefore we obtain:

\begin{corollary}
Assume that $(X_{k})$ is as in Theorem \ref{Th attraction}, (\ref{cond0}),
(\ref{tp}) and (\ref{lrtp}) are satisfied. Then the conclusions of Theorem
\ref{Th attraction} hold.
\end{corollary}

\subsection{Examples}

\textbf{Example 1. }The following is an example of a strictly stationary
Markov chain satisfying (\ref{main condition})\ but not the right hand side of
(\ref{cond0}) (see Bradley (1997), Remark 1.5). The Markov chain $(\xi
_{k})_{k\in{\mathbb{Z}}}$ is such that the marginal distribution of $\xi_{0}$
is uniformly distributed on $[0,1]$ and the one step transition probabilities
are as follows:

For each $x\in\lbrack0,1],$ let
\[
\mathbb{P}(\xi_{1}=x|\xi_{0}=x)=1/2
\]
and%
\[
\mathbb{P}(\xi_{1}\in B|\xi_{0}=x)=(1/2)\lambda(B),
\]
where $B\subset\lbrack0,1]-x$ and $\lambda$ is the Lebesgue measure on
$[0,1]$. Then (\ref{main condition}) holds with $\ a=1/2.$ However no such $b$
as in the right hand side of (\ref{cond0}) exists.

\bigskip

\textbf{Example 2. (}Generalized\textbf{ }Gibbs Markov chains\textbf{). }Let
$S$ be a countable set, $p:S\times S\rightarrow\lbrack0,1]$ be an aperiodic,
irreducible stochastic matrix and $(\pi_{s})_{s\in S}$ an invariant
distribution with $\pi_{s}>0$ for all $s\in S$. Let $T:S^{{\mathbb{N}}%
}\rightarrow S^{{\mathbb{N}}}$ be the shift and define the Markov chain in a
canonical way on $S^{{\mathbb{N}}}$ by%
\[
{\mathbb{P}}(X_{1}=x_{1},...,X_{n}=x_{n})=\pi_{x_{1}}p(x_{1},x_{2}%
)...p(x_{n-1},x_{n}).
\]
Let $\Omega=\{x\in S^{{\mathbb{N}}}:{\mathbb{P}}(X_{1}=x_{1},...,X_{n}%
=x_{n})>0\}.$ We assume that there is $0<m<1$ such that for all $s,t\in S$
\begin{equation}
p(s,t)\geq m\pi_{t} \label{Gibbs}%
\end{equation}
Then our condition (\ref{main condition}) is satisfied with $a=m$. This is a
larger class than the so called Gibbs-Markov maps as defined in Aaronson and
Denker (2001b).

\bigskip

\textbf{Example 3}. In the context of Example 2, a fairly large class of
countable state Markov processes satisfying condition (\ref{Gibbs}) can be
constructed by defining for $i,j\in{\mathbb{N}}^{\ast}$%
\[
p(i,j)=\pi_{j}+(\delta_{i,j}-\delta_{i+1,j})\varepsilon_{i},
\]
where for all $i\in{\mathbb{N}}^{\ast},$ $\delta_{i,i}=1$ and for $j\neq i$ we
have $\delta_{i,j}=0.$ We take $0\leq\varepsilon_{i}\leq\min(1-\pi_{i}%
,\pi_{i+1}).$ In addition we assume that there is $0<m<1$ such that
\[
\varepsilon_{i}\leq(1-m)\pi_{i+1}.
\]
For example, let $m=1/2,$ $\pi_{j}=2^{-j}$ and set $p(i,j)=2^{-j}%
+(\delta_{i,j}-\delta_{i+1,j})2^{-(2+i)}.$

\section{Relation to mixing coefficients\label{mix}}

\textbf{Relation to lower }$\psi-$\textbf{mixing coefficient.}

We introduce now a mixing condition which is comparable to condition
(\ref{main condition}).

Following Bradley (2007), for any two sigma algebras $\mathcal{A}$ and
$\mathcal{B}$ define the lower $\psi-$mixing coefficient by
\[
\psi^{\prime}(\mathcal{A},\mathcal{B})=\inf\frac{\mathbb{P}(A\cap
B)}{\mathbb{P}(A)\mathbb{P}(B)};\text{ }A\in\mathcal{A}\text{ and }%
B\in\mathcal{B}\text{, }\mathbb{P}(A)\mathbb{P}(B)>0.
\]
Obviously $0\leq\psi^{\prime}(\mathcal{A},\mathcal{B})\leq1$ and $\psi
^{\prime}(\mathcal{A},\mathcal{B})=1$ if and only if $\mathcal{A}$ and
$\mathcal{B}$ are independent.

For a sequence $(\xi_{k})_{k\geq0}$ of random variables we shall denote by
$\mathcal{F}_{n}^{m}=\sigma(\xi_{i},m\leq i\leq n)$ and by
\[
\psi_{k}^{\prime}=\inf_{m\geq0}\psi^{\prime}(\mathcal{F}_{0}^{m}%
,\mathcal{F}_{k+m}^{k+m}).
\]
In the Markov setting, by the Markov property,
\[
\psi_{k}^{\prime}=\inf_{m\geq0}\psi^{\prime}(\sigma(\xi_{m}),\sigma(\xi
_{k+m})).
\]
Notice that, in terms of conditional probabilities defined by (\ref{defQ}), we
also have the following equivalent definition:
\[
\psi_{1}^{\prime}=\inf_{k\geq1}\mathrm{ess}\inf_{x}\inf_{A\in\mathcal{B}%
(\mathcal{X})}Q_{k}(x,A)/\mathbb{P}(\xi_{k}\in A)\text{.}%
\]
Note that if (\ref{main condition})\ holds then $\psi_{1}^{\prime}\geq a>0$
and if $\psi_{1}^{\prime}>0$ then (\ref{main condition})\ holds with
$a=\psi_{1}^{\prime}.$ Also note that condition (\ref{main condition}) is
equivalent to the existence of a constant $a^{\prime}$ such that $\psi
_{1}^{\prime}=a^{\prime}>0.$

For Markov chains, by Theorem 7.4 (d) in Bradley (2007):%
\begin{equation}
1-\psi_{k+m}^{\prime}\leq(1-\psi_{k}^{\prime})(1-\psi_{m}^{\prime}).
\label{prod}%
\end{equation}
So, by condition (\ref{main condition}), we have $1-\psi_{n}^{\prime}%
\leq(1-a^{\prime})^{n}\rightarrow0,$ and in this case, $(\xi_{k})_{k\geq0}$ is
called lower $\psi$-mixing.

It should be mentioned that the condition on the coefficient $\psi_{1}%
^{\prime}\ $is quite natural for obtaining large deviation results for
stationary Markov chains (Bryc and Smolenski, 1993).

\bigskip

\textbf{Relation to }$\rho-$\textbf{mixing coefficients.}

\bigskip

Define the maximal coefficient of correlation
\[
\rho(\mathcal{A},\mathcal{B})=\sup_{X\in{\mathbb{L}}_{2}(\mathcal{A}%
),Y\in{\mathbb{L}}_{2}(\mathcal{B})}|\mathrm{corr}\,(X,Y)|\text{ .}%
\]

By Lemma 10 in Merlev\`{e}de et al. (2021) (which is actually due to R.
Bradley), we know that%
\begin{equation}
\rho(\mathcal{A},\mathcal{B})\leq1-\psi^{\prime}(\mathcal{A},\mathcal{B}).
\label{key}%
\end{equation}
For a sequence $(\xi_{k})_{k\geq0}$ of random variables $\rho_{k}%
\mathbf{=}\sup_{m\geq0}\rho(\mathcal{F}_{0}^{m},\mathcal{F}_{k+m}^{\infty}),$
and in the Markov case $\rho_{k}\mathbf{=}\sup_{m\geq0}\rho(\sigma(\xi
_{m}),\sigma(\xi_{k+m})).$ In this latter case we also have the equivalent
definition for $\rho_{1},$
\[
\rho_{_{1}=}\sup_{k\geq1}\sup_{||f||_{2}=1,\mathbb{E}f=0}||Q_{k}f||_{2}\text{
},
\]
where%
\[
Q_{k}f(x)=\int f(y)Q_{k}(x,dy).
\]
According to (\ref{key}), condition (\ref{main condition})\ implies
\begin{equation}
\rho_{1}\leq1-a<1. \label{relro}%
\end{equation}
By Theorem 7.4 (a) in Bradley (2007), it follows that $\rho_{n}\leq
(1-a)^{n}\rightarrow0.$ Therefore, a Markov chain satisfying condition
(\ref{main condition}) is also $\rho$-mixing.

Define now $(X_{k})_{k\geq0}$ by (\ref{defX}). Assume the variables
$(X_{k})_{k\geq0}$ are centered and have finite second moments. Denote
$\tau_{n}^{2}=\sum_{j=1}^{n}\mathrm{var}(X_{j})$ and $\sigma_{n}%
^{2}={\mathbb{E}}(S_{n}^{2})$. From Proposition 13 in Peligrad (2012) we know
that
\[
\frac{1-\rho_{1}}{1+\rho_{1}}\leq\frac{\sigma_{n}^{2}}{\tau_{n}^{2}}\leq
\frac{1+\rho_{1}}{1-\rho_{1}}.
\]
By combining this inequality with (\ref{relro}), we obtain
\begin{equation}
\frac{a}{2-a}\leq\frac{\sigma_{n}^{2}}{\tau_{n}^{2}}\leq\frac{2-a}{a}.
\label{varineq}%
\end{equation}
If the Markov chain is stationary define $X_{k}=g(\xi_{k})$ and assume that
the variables $X_{k}$ are nondegenerate, have finite second moment and are
centered. Then, by the definition of $\rho_{k}$ and (\ref{relro}), for all
$k>0,$ $|\mathbb{E}(X_{0}X_{k})|\leq(1-a)^{k}\mathbb{E}(X_{k}^{2})$\ and
therefore\ it follows that there is $c>0$ such that%
\begin{equation}
\frac{\sigma_{n}^{2}}{n}\rightarrow c^{2}. \label{limvar}%
\end{equation}

\textbf{Relation to} $\varphi-$\textbf{mixing coefficients}

\bigskip

Relevant to our paper is also the relation with $\varphi-$mixing coefficient
defined by
\[
\varphi(\mathcal{A},\mathcal{B})=\sup\frac{\mathbb{P}(A\cap B)-\mathbb{P}%
(A)\mathbb{P}(B)}{\mathbb{P}(A)};\text{ }A\in\mathcal{A}\text{ and }%
B\in\mathcal{B}\text{, }\mathbb{P}(A)>0.
\]
By Proposition 5.2 (III)(b) in Bradley (2007),
\begin{equation}
\varphi(\mathcal{A},\mathcal{B})\leq1-\psi^{\prime}(\mathcal{A},\mathcal{B}).
\label{compphy}%
\end{equation}
In the Markov case
\[
\varphi_{k}=\sup_{m\geq0}\varphi(\sigma(\xi_{m}),\sigma(\xi_{k+m})).
\]
Then, if condition (\ref{main condition}) is satisfied, by (\ref{compphy}) and
by (\ref{prod}), for any $k\geq1$
\begin{equation}
\varphi_{k}\leq1-\psi_{k}^{\prime}\leq(1-a)^{k}\rightarrow0\text{ as
}n\rightarrow\infty, \label{phi1sm1}%
\end{equation}
and so, our Markov chain is also $\varphi-$mixing, at least exponentially fast.

\bigskip

\textbf{Relation to }$\psi-$\textbf{mixing coefficient.}

\bigskip

The class of $\psi^{\prime}-$mixing Markov chains is strictly larger that the
class of $\psi-$mixing. To see this we mention that the $\psi-$mixing
coefficient is defined in the following way. For any two sigma algebras
$\mathcal{A}$ and $\mathcal{B}$, define the $\psi-$mixing coefficient by
\[
\psi(\mathcal{A},\mathcal{B})=\sup\frac{\mathbb{P}(A\cap B)-\mathbb{P}%
(A)\mathbb{P}(B)}{\mathbb{P}(A)\mathbb{P}(B)};\text{ }A\in\mathcal{A}\text{
and }B\in\mathcal{B}\text{, }\mathbb{P}(A)\mathbb{P}(B)>0.
\]
It is well-known that%
\[
\psi(\mathcal{A},\mathcal{B})=\max[\psi^{\ast}(\mathcal{A},\mathcal{B}%
)-1,1-\psi^{\prime}(\mathcal{A},\mathcal{B})],
\]
where%
\[
\psi^{\ast}=\sup\frac{\mathbb{P}(A\cap B)}{\mathbb{P}(A)\mathbb{P}(B)};\text{
}A\in\mathcal{A}\text{ and }B\in\mathcal{B}\text{, }\mathbb{P}(A)\mathbb{P}%
(B)>0.
\]
For a sequence $(\xi_{k})_{k\geq0}$ of random variables $\psi_{k}=\sup
_{m\geq0}\psi(\mathcal{F}_{0}^{m},\mathcal{F}_{k+m}^{k+m})\ $and in the
Markovian case $\psi_{k}=\sup_{m\geq0}\psi(\sigma(\xi_{m}),\sigma(\xi
_{k+m})).$ Note that $\psi_{1}<1$ implies $\psi_{1}^{\prime}>0$. Therefore all
our results also hold for $\psi-$mixing Markov chains with $\psi_{1}<1.$

Furthermore, the class of $\psi$-mixing Markov chains is strictly smaller than
the class of $\psi^{\prime}$-mixing Markov chains (see for example Remark 1.5
in Bradley, 1997, and Example 1 above). Therefore, clearly, the conclusions of
our results also hold for the smaller class of $\psi$-mixing Markov chains
with $\psi_{1}<1$.

Note that condition (\ref{cond0}) is equivalent to the existence of two
constants $a^{\prime}$ and $b^{\prime}$ such that $\psi_{1}^{\prime}%
=a^{\prime}>0$ and $\psi_{1}^{\ast}=b^{\prime}<\infty.$

\section{Proofs}

\subsection{Preliminary general local CLT}

Here we give two general local limit theorems.

\begin{theorem}
\label{general}We require that convergence in distribution in
(\ref{conv distribution}) holds. In addition, suppose that for each $D>0$
\begin{equation}
\ \lim_{T\rightarrow\infty}\text{\ }\limsup_{n\rightarrow\infty}%
\int\nolimits_{T<|t|\leq DB_{n}}\left\vert {\mathbb{E}}\exp\left(
it\frac{S_{n}}{B_{n}}\right)  \right\vert dt=0\text{.} \label{integral}%
\end{equation}
Then (\ref{GLLTnonlat}) holds.
\end{theorem}

\begin{remark}
\label{remd1d2}By decomposing the integral in (\ref{integral}) into two parts,
on $\{T<|u|<\delta B_{n}\}$ and on $\{\delta B_{n}\leq|u|\leq DB_{n}\}$, and
changing the variable in the second integral we easily argue that in order to
prove this theorem it is enough to show that for each $D$ fixed there is
$0<\delta<D$ such that%
\[
(D_{1})\text{ \ \ \ }\lim_{T\rightarrow\infty}\limsup_{n\rightarrow\infty}%
\int\nolimits_{T<|t|<\delta B_{n}}\left\vert {\mathbb{E}}\exp\left(
it\frac{S_{n}}{B_{n}}\right)  \right\vert dt=0
\]
and $\ $%
\[
(D_{2})\text{ \ \ \ \ \ \ \ \ \ \ \ \ \ \ }\lim_{n\rightarrow\infty}B_{n}%
\int\nolimits_{\delta\leq|t|\leq D}|{\mathbb{E}}\exp(itS_{n})|dt=0.
\]
\ \ \ \ \ 
\end{remark}

If all the variables $(S_{n})$ have the values in a fixed lattice, we impose a
different condition and the result we obtain is the following:

\begin{theorem}
\label{ThLattice}Assume now that all the variables $(S_{n})$ have values in a
fixed lattice%
\[
\mathcal{S}=\{kh,k \in{\mathbb{Z}}\}
\]
and (\ref{conv distribution}) holds. Suppose in addition that
\begin{equation}
\lim_{T\rightarrow\infty}\limsup_{n\rightarrow\infty}\int_{T<|t|\leq\frac{\pi
}{h}B_{n}}\left\vert {\mathbb{E}}\exp\left(  it\frac{S_{n}}{B_{n}}\right)
\right\vert dt=0. \label{Mcond}%
\end{equation}
Then, the conclusion in (\ref{GLLTlatt}) holds.
\end{theorem}

\begin{remark}
By decomposing the integral in (\ref{Mcond}) into two parts, on
$\{T<|u|<\delta B_{n}\}$ and on $\{\delta B_{n}\leq|u|\leq\pi B_{n}/h\}$, the
conclusion of Theorem \ref{ThLattice} holds if, for some $0<\delta<\pi/h$,
\[
(D_{1}^{\prime})\text{ \ \ \ }\lim_{T\rightarrow\infty}\limsup_{n\rightarrow
\infty}\int\nolimits_{T<|t|<\delta B_{n}}\left\vert {\mathbb{E}}\exp\left(
it\frac{S_{n}}{B_{n}}\right)  \right\vert dt=0
\]
and
\[
(D_{2}^{\prime})\text{ \ \ \ \ \ \ \ }\lim_{n\rightarrow\infty}B_{n}%
\int\nolimits_{\delta\leq|t|\leq\pi/h}|{\mathbb{E}}\exp(itS_{n})|dt=0.\text{
\ \ \ \ \ }%
\]

\end{remark}

\textbf{Proof of Theorems \ref{general} and \ref{ThLattice}.} The proofs of
these theorems can be deduced from the corresponding arguments in Hafouta and
Kifer (2016), who assumed convergence to normal distribution and normalization
$\sqrt{n}.$ Since the proof of Theorem \ref{general} is given in Peligrad et
al. (2021) we shall give here only the proof of Theorem \ref{ThLattice}.

The proof is based on the inverse Fourier formula for the sum $S_{n}$.
According to Lemma 4.5 and arguments in Section VI.4 in Hennion and Herv\'{e}
(2001) (see also Theorem 10.7 in Breiman (1982) and Section 10.4 there), it
suffices to prove (\ref{GLLTlatt}) for all continuous complex valued functions
$g$ defined on ${\mathbb{R}}$ such that $|g|\in L^{1}({\mathbb{R}})$ and its
Fourier transform
\[
\hat{g}(t)=\int\nolimits_{{\mathbb{R}}}e^{-itx}g(x)\ dx
\]
has compact support contained in some finite interval $[-D,D]$.

The inversion formula gives:%
\[
g(x)=\frac{1}{2\pi}\int\nolimits_{{\mathbb{R}}}e^{isx}\hat{g}(s)\ ds.
\]
Therefore, using a change of variable and taking the expected value, we obtain%

\[
{{\mathbb{E}}\,}[g(S_{n}+u)]=\frac{1}{2\pi B_{n}}\int_{{\mathbb{R}}}\hat
{g}\Big (\frac{t}{B_{n}}\Big )\ f_{S_{n}}\Big (\frac{t}{B_{n}}\Big )\ \exp
\Big (\frac{itu}{B_{n}}\Big )\ dt,
\]
where we have used the notation%
\[
f_{X}(v)=\mathbb{E}(\exp(ivX)),
\]
and $u\in\mathcal{S}$, where $\mathcal{S}$ is as in the statement of Theorem
\ref{ThLattice}. Denote
\[
r(v)=\sum_{k=-\infty}^{\infty}\hat{g}\Big (v+2\pi\frac{k}{h}\Big ).
\]
Since $\hat{g}$ has compact support included in an interval $[-D,D]$, the sum
is finite. So, with the notation $M_{v}=h(D+|v|)/2\pi$, clearly%
\begin{equation}
|r(v)|\leq\sum_{|k|\leq M_{v}}\Big |\hat{g}\Big (v+2\pi\frac{k}{h}%
\Big )\Big |\leq2M_{v}\sup_{|w|\leq D}|\hat{g}(w)|.\label{prop r}%
\end{equation}
Now,%
\begin{align*}
&  \int_{{\mathbb{R}}}\hat{g}\Big(\frac{t}{B_{n}}\Big)\;f_{S_{n}}\Big(\frac
{t}{B_{n}}\Big)\;\exp\Big(\frac{itu}{B_{n}}\Big)\ dt\\
&  =\sum_{k=-\infty}^{\infty}\int_{\frac{2k-1}{h}\pi B_{n}}^{\frac{2k+1}{h}\pi
B_{n}}\hat{g}\Big(\frac{t}{B_{n}}\Big)\;f_{S_{n}}\Big(\frac{t}{B_{n}%
}\Big)\;\exp\Big(\frac{itu}{B_{n}}\Big)\ dt.
\end{align*}
We change the variable from $\ t$ to $v$ according to:%
\[
\frac{t}{B_{n}}=\frac{v}{B_{n}}+\frac{2\pi k}{h}.
\]
Because $u=k^{\prime}h,$ $k^{\prime}\in{\mathbb{Z}}$,
\[
\exp\Big (\frac{itu}{B_{n}}\Big )=\exp\Big  (iu\frac{v}{B_{n}}+ik^{\prime
}h\frac{2\pi k}{h}\Big )=\exp\Big  (iu\frac{v}{B_{n}}\Big  ).
\]
Since $S_{n}$ has values in $\mathcal{S}$ we note that
\[
f_{S_{n}}\Big (\frac{t}{B_{n}}+\frac{2k}{h}\pi\Big )=f_{S_{n}}\Big (\frac
{t}{B_{n}}\Big )\text{ for all }k\in{\mathbb{Z}}.
\]
Therefore, by using the definition of $r(v),$%
\[
{{\mathbb{E}}\,}[g(S_{n}+u)]=\frac{1}{2\pi B_{n}}\int_{-\frac{\pi}{h}B_{n}%
}^{\frac{\pi}{h}B_{n}}f_{S_{n}}\Big  (\frac{v}{B_{n}}\Big )\;\exp
\Big  (\frac{ivu}{B_{n}}\Big )\ r\Big (\frac{v}{B_{n}}\Big )\ dv.
\]
Note that, by the definition of the characteristic function, we have the
identity $f_{L}(u)=\hat{h}_{L}(-u)$. So, by the Fourier inversion formula we
also have%
\[
h_{L}\Big  (-\frac{u}{B_{n}}\Big  )=\frac{1}{2\pi}\int f_{L}(t)\text{ }%
\exp\Big (\frac{itu}{B_{n}}\Big )\ dt.
\]
\ We evaluate
\begin{align*}
&  2\pi\left\vert B_{n}\mathbb{E}[g(S_{n}+u)]-h_{L}(-\frac{u}{B_{n}}%
)\sum_{k=-\infty}^{\infty}hg(kh)\right\vert \\
&  =\left\vert \int_{-\frac{\pi}{h}B_{n}}^{\frac{\pi}{h}B_{n}}f_{S_{n}%
}\Big (\frac{v}{B_{n}}\Big )\;\exp\Big (\frac{ivu}{B_{n}}\Big )\ r\Big (\frac
{v}{B_{n}}\Big )\ dv.-\int f_{L}(v)\text{ }\exp\Big (\frac{ivu}{B_{n}%
}\Big )\ dv\int gd\mathcal{L}_{h}\right\vert \\
&  =\left\vert I_{n}+II_{n}+III_{n}\right\vert \leq|I_{n}|+|II_{n}|+|III_{n}|,
\end{align*}
where the terms $I_{n},$ $II_{n},$ and $III_{n}$ are given below and their
moduli analyzed:%
\begin{align*}
|I_{n}| &  =\left\vert \int_{-T}^{T}\left(  f_{S_{n}}\Big (\frac{t}{B_{n}%
}\Big )\ r\Big (\frac{t}{B_{n}}\Big )-f_{L}(t)\int gd\mathcal{L}_{h}\right)
\text{ }\exp(\frac{itu}{B_{n}})\ dt\right\vert \\
&  \leq\int_{-T}^{T}\left\vert f_{S_{n}}\Big (\frac{t}{B_{n}}%
\Big )r\Big (\frac{t}{B_{n}}\Big )\ dv-f_{L}(t)\int gd\mathcal{L}%
_{h}\right\vert dt,
\end{align*}%
\[
|II_{n}|\leq\left(  \int_{|t|>T}|f_{L}(t)|dt\right)  \left(  |\int
gd\mathcal{L}_{h}|\right)  ,
\]
and%
\[
|III_{n}|\leq\int_{T<|v|\leq\frac{\pi}{h}B_{n}}\Big |f_{S_{n}}\Big(\frac
{v}{B_{n}}\Big )\Big |\ \Big |r\Big (\frac{v}{B_{n}}\Big )\Big |\ dv.
\]
To deal with $|I_{n}|,$ we further decompose the sum and use the triangle
inequality,%
\begin{align*}
|I_{n}| &  \leq\int_{-T}^{T}\left\vert f_{S_{n}}\Big(\frac{v}{B_{n}%
}\Big )r\Big (\frac{v}{B_{n}}\Big)-f_{S_{n}}\Big (\frac{v}{B_{n}}\Big )\int
gd\mathcal{L}_{h}\right\vert \ dv\\
&  +\int_{-T}^{T}\left\vert f_{S_{n}}\Big (\frac{v}{B_{n}}\Big )\int
gd\mathcal{L}_{h}-f_{L}(v)\int gd\mathcal{L}_{h}\right\vert dv.
\end{align*}
The first term in the right hand side of the above inequality tends to $0$
because by Ch 10 in Breiman (1992) and since $B_{n}\rightarrow\infty,$ we have%
\[
r\Big(\frac{v}{B_{n}}\Big)\rightarrow\int gd\mathcal{L}_{h}.
\]
The second term in the right hand side converges to $0$ because, by our
hypotheses
\[
f_{S_{n}}\Big (\frac{v}{B_{n}}\Big )\rightarrow f_{L}(v),
\]
uniformly on compacts.

The term $|II_{n}|$ converges to $0$ as $T\rightarrow\infty$, because $f_{L}$
is assumed to be integrable.

Next, for $|III_{n}|$ we have
\begin{align*}
|III_{n}| &  \leq\sup_{T<|w|\leq\frac{\pi}{h}B_{n}}\Big |r\Big (\frac{w}%
{B_{n}}\Big )\Big |\int_{T<|v|\leq\frac{\pi}{h}B_{n}}\Big |f_{S_{n}%
}\Big (\frac{v}{B_{n}}\Big )\Big |dv\\
&  \leq\sup_{|w|\leq\frac{\pi}{h}}|r(w)|\int_{T<|v|\leq\frac{\pi}{h}B_{n}%
}\Big |f_{S_{n}}\Big (\frac{v}{B_{n}}\Big )\Big |dv.
\end{align*}
The result follows by condition (\ref{Mcond}) and the fact that by
(\ref{prop r}) we have
\[
\sup_{|w|\leq\frac{\pi}{h}}|r(w)|\leq\left[  (Dh+\pi)/\pi\right]
\sup_{|w|\leq D}|\hat{g}(w)|.
\]

$\square$

\subsection{Bounds on the characteristic function}

The bound on the characteristic function of a Markov chain is inspired by
Lemma 1.5 in Nagaev (1961). It is given in the following proposition. Recall
the definition of $a$ in (\ref{main condition}).

\begin{proposition}
\label{factor}Let $(X_{j})_{j\geq0}$ be defined by (\ref{defX}). Then, for all
$n\in\mathbb{N},$ and $u\in\mathbb{R}$,
\[
|{\mathbb{E}}(\exp(iuS_{n}))|\leq\exp\Big [-\frac{a^{4}}{16}\sum_{j=1}%
^{n}(1-|f_{j}(u)|^{2})\Big ].
\]

\end{proposition}

For proving this proposition we need some preliminary considerations.

For a complex, finite measure $\mu,$ defined on a sigma algebra $\mathcal{F}$,
denote by \textrm{Var}$\mu$ its total variation. This is, for all measurable
$A$ we have%

\[
\mathrm{Var}\mu(A)=\sup\sum\nolimits_{k=1}^{n}|\mu(A_{k})|,
\]
where the supremum is taken over all $n\in\mathbb{N},$ and all $A_{1},$
$A_{2},\ldots,A_{n,}$ disjoint sets in $\mathcal{F}$, with $A=\cup_{k=1}%
^{n}A_{k}.$

The property we shall use below is that for complex measurable functions
integrable with respect to $\mu$ we have%
\[
\Big |\int h(y)d\mu\Big |\leq\int|h(y)|d(\mathrm{Var}\mu).
\]

Denote by $\mathbb{L}_{\infty}(\mathcal{X},\mathcal{B}(\mathcal{X}%
),\mathbb{P}^{\prime})$ the space of complex valued, measurable, essentially
bounded functions on the probability space $(\mathcal{X},\mathcal{B}%
(\mathcal{X}),\mathbb{P}^{\prime}).$ For an operator $T:$ $\mathbb{L}_{\infty
}(\mathcal{X},\mathcal{B}(\mathcal{X}),\mathbb{P}^{\prime})\rightarrow
\mathbb{L}_{\infty}(\mathcal{X},\mathcal{B}(\mathcal{X}),\mathbb{P}^{"})$ we
denote by%
\[
||T||:=\sup_{||f||_{1}=1}||T(f)||_{1},
\]
where $||T(f)||_{1}$ is the norm of $T(f)$ in the space of complex valued,
integrable functions on $(\mathcal{X},\mathcal{B}(\mathcal{X}),\mathbb{P}")$
and $||f||_{1}$ is the norm of $f$ in the space of complex valued integrable
functions on $(\mathcal{X},\mathcal{B}(\mathcal{X}),\mathbb{P}^{{\prime}}%
)$.\bigskip

We introduce now a family of operators which are relevant to our proofs:

For $u$ a fixed real number let us introduce the operator
\[
T_{u,k}:\mathbb{L}_{\infty}(\mathcal{X},\mathcal{B}(\mathcal{X}),\mathbb{P}%
_{k})\rightarrow\mathbb{L}_{\infty}(\mathcal{X},\mathcal{B}(\mathcal{X}%
),\mathbb{P}_{k-1})
\]
by
\[
T_{u,k}(h)(x):=\int h(y)\exp(iug_{k}(y))Q_{k}(x,dy).
\]
So, for $k\geq1$,%
\[
T_{u,k}(h)(\xi_{k-1})={\mathbb{E}}\left(  [h(\xi_{k})\exp(iuX_{k}%
)]{\LARGE |}\xi_{k-1}\right)  .
\]

\begin{lemma}
\label{estimate3} For any $k\in{\mathbb{N}}^{*}$ and $u\in\mathbb{R}$ we
have,
\[
||T_{u,k-1}\circ T_{u,k}||_{1}\leq1-\frac{a^{4}}{8}(1-|{\mathbb{E}}%
(\exp(iuX_{k-1})|^{2}).
\]

\end{lemma}

\noindent\textbf{Proof.} In the sequel, to simplify the notation we drop the
index $u,$ and write $T_{k}=T_{u,k}.$ Let $x\in\mathcal{X}^{\prime},$ where
$\mathcal{X}^{\prime}\in$ $\mathcal{B}(\mathcal{X})$ such that $\mathbb{P}%
_{k-1}(\mathcal{X}^{\prime})=1,$ for which condition (\ref{cond0}) holds. By
the definition of $T_{k}$'s,
\begin{align*}
&  T_{k-1}\circ T_{k}(h)(x)\\
&  =\int\exp(iug_{k-1}(y))\int h(z)\exp(iug_{k}(z))Q_{k}(y,dz)Q_{k-1}(x,dy).
\end{align*}
Changing the order of integration
\begin{align*}
&  T_{k-1}\circ T_{k}(h)(x)\\
&  =\int h(z)\exp(iug_{k}(z))\int\exp(iug_{k-1}(y))Q_{k-1}(x,dy)Q_{k}(y,dz)\\
&  =\int h(z)\exp(iug_{k}(z))m_{x}(dz),
\end{align*}
where, for $x$ fixed in $\mathcal{X}^{\prime}$, $m_{x}$ is the complex finite
measure defined on $\mathcal{B}(\mathcal{X})$ by the formula
\[
m_{x}(A)=\int\exp(iug_{k-1}(y))Q_{k}(y,A)Q_{k-1}(x,dy).
\]
Note that%
\begin{multline}
\int|T_{k-1}\circ T_{k}(h)(x)|\mathbb{P}_{k-2}(dx)=\int\left\vert \int
h(z)\exp(iug_{k}(z))m_{x}(dz)\right\vert \mathbb{P}_{k-2}(dx)\label{rep T}\\
\leq\int\int|h(z)|\mathrm{Var}\left(  m_{x}(dz)\right)  \mathbb{P}%
_{k-2}(dx):=\int|h(z)|m(dz),
\end{multline}
where
\[
m(A)=\int\mathrm{Var}\left(  m_{x}(A)\right)  \mathbb{P}_{k-2}(dx).
\]
To analyze $m(A)$, we start by noting that for an increasing sequence of
measurable partitions $\mathcal{P}_{n}$ of $A,$ by the monotone convergence
theorem,
\begin{equation}
m(A)=\lim_{n\rightarrow\infty}\sum\nolimits_{A_{i}\in\mathcal{P}_{n}}%
\int|m_{x}(A_{i})|\mathbb{P}_{k-2}(dx).\label{def m(A)}%
\end{equation}
In order to find an upper bound for $|m_{x}(A)|$, we start from the following
estimate, where we have used at the end condition (\ref{main condition}):%
\begin{gather*}
\left(  \int Q_{k}(y,A)Q_{k-1}(x,dy)\right)  ^{2}-\left(  \Big |\int
\exp(iug_{k-1}(y))Q_{k}(y,A)Q_{k-1}(x,dy)\Big |\right)  ^{2}\\
=\iint(1-\cos(u\left(  g_{k-1}(y)-g_{k-1}(y^{\prime}))\right)  Q_{k}%
(y,A)Q_{k-1}(x,dy)Q_{k}(y^{\prime},A)Q_{k-1}(x,dy^{\prime})\\
=\iint2\sin^{2}\left(  \frac{u}{2}(g_{k-1}(y)-g_{k-1}(y^{\prime}))\right)
Q_{k}(y,A)Q_{k-1}(x,dy)Q_{k}(y^{\prime},A)Q_{k-1}(x,dy^{\prime})\\
\geq a^{4}\mathbb{P}_{k}^{2}(A)M,
\end{gather*}
where%
\[
M=2\iint\sin^{2}\left(  \frac{u}{2}(g_{k-1}(y)-g_{k-1}(y^{\prime}))\right)
\mathbb{P}_{k-1}(dy)\mathbb{P}_{k-1}(dy^{\prime}).
\]
Assume for the moment that $\mathbb{P}_{k}(A)>0.$ By condition
(\ref{main condition}) this implies that, for all $x\in\mathcal{X}^{\prime}$
\[
\int Q_{k}(y,A)Q_{k-1}(x,dy)\geq a\mathbb{P}_{k}(A)>0.
\]
Because%
\begin{align*}
|m_{x}(A)| &  =\int Q_{k}(y,A)Q_{k-1}(x,dy)-\int Q_{k}(y,A)Q_{k-1}(x,dy)\\
&  +\Big |\int\exp(iug_{k-1}(y))Q_{k}(y,A)Q_{k-1}(x,dy)\Big |,
\end{align*}
we have%
\begin{gather*}
|m_{x}(A)|=\int Q_{k}(y,A)Q_{k-1}(x,dy)\\
-\frac{\left(  \int Q_{k}(y,A)Q_{k-1}(x,dy)\right)  ^{2}-\left(  |\int
\exp(iug_{k-1}(y))Q_{k}(y,A)Q_{k-1}(x,dy)|\right)  ^{2}}{\int Q_{k}%
(y,A)Q_{k-1}(x,dy)+|\int\exp(iug_{k-1}(y))Q_{k}(y,A)Q_{k-1}(x,dy)|}.
\end{gather*}
But%
\begin{align*}
&  \int Q_{k}(y,A)Q_{k-1}(x,dy)+\Big |\int\exp(iug_{k-1}(y))Q_{k}%
(y,A)Q_{k-1}(x,dy)\Big |\\
&  \leq2\int Q_{k}(y,A)Q_{k-1}(x,dy).
\end{align*}
Therefore, by the above considerations,
\[
|m_{x}(A)|\leq\int Q_{k}(y,A)Q_{k-1}(x,dy)-\frac{a^{4}\mathbb{P}_{k}^{2}%
(A)M}{2\int Q_{k}(y,A)Q_{k-1}(x,dy)}.
\]
With the notation
\[
Y_{k-2}(x)=\int Q_{k}(y,A)Q_{k-1}(x,dy)=\mathbb{P}(\xi_{k}\in A|\xi_{k-2}=x),
\]
by integrating the above inequality with respect to $\mathbb{P}_{k-2}(dx)$ we
obtain
\begin{gather}
\int|m_{x}(A)|\mathbb{P}_{k-2}(dx)\leq\int Y_{k-2}(x)\mathbb{P}_{k-2}%
(dx)\label{ineq int m(x,A)}\\
-a^{4}\mathbb{P}_{k}(A)M\int\frac{\mathbb{P}_{k}(A)}{2Y_{k-2}(x)}%
\mathbb{P}_{k-2}(dx).\nonumber
\end{gather}
By the Markov inequality,
\[
\mathbb{P}_{k-2}(Y_{k-2}>2\mathbb{P}_{k}(A))\leq\frac{\int Y_{k-2}%
(x)\mathbb{P}_{k-2}(dx)}{2\mathbb{P}_{k}(A)}=\frac{\mathbb{P}_{k}%
(A)}{2\mathbb{P}_{k}(A)}=\frac{1}{2}.
\]
So%
\[
\mathbb{P}_{k-2}(Y_{k-2}\leq2\mathbb{P}_{k}(A))\geq\frac{1}{2}.
\]
Therefore%
\begin{multline*}
\int\frac{\mathbb{P}_{k}(A)}{Y_{k-2}(x)}\mathbb{P}_{k-2}(dx)\geq\int
_{Y_{k-2}\leq2\mathbb{P}_{k}(A)}\frac{\mathbb{P}_{k}(A)}{Y_{k-2}(x)}%
\mathbb{P}_{k-2}(dx)\\
\geq\frac{1}{2}\mathbb{P}_{k-2}(Y_{k-2}\leq2\mathbb{P}_{k}(A))\geq\frac{1}{4}.
\end{multline*}
Combining the last two inequalities with (\ref{ineq int m(x,A)}) we obtain
\begin{equation}
\int|m_{x}(A)|\mathbb{P}_{k-2}(dx)\leq\mathbb{P}_{k}(A)-\frac{1}{8}%
a^{4}\mathbb{P}_{k}(A)M=\left(  1-\frac{1}{8}a^{4}M\right)  \mathbb{P}%
_{k}(A).\label{ineq mxa}%
\end{equation}
Let us note now that inequality (\ref{ineq mxa}) is also true for
$\mathbb{P}_{k}(A)=0.$ Indeed, since
\[
\mathbb{P}_{k}(A)=\int Q_{k}(y,A)\mathbb{P}_{k-1}(dy),
\]
we get $Q_{k}(y,A)=0,$ $\mathbb{P}_{k-1}-a.s.$, and by its definition,
$m_{x}(A)=0$ for all $x\in\mathcal{X}^{\prime}$.

At this moment, inequality (\ref{ineq mxa})\ allows to estimate $m(A)$ defined
in (\ref{def m(A)}). Hence, we get
\begin{multline*}
m(A)=\lim_{n}\ \sum\nolimits_{A_{i}\in\mathcal{P}_{n}}\int|m_{x}%
(A_{i})|\mathbb{P}_{k-2}(dx)\\
\leq\lim_{n}\ \sum\nolimits_{A_{i}\in\mathcal{P}_{n}}\left(  1-\frac{1}%
{8}a^{4}M\right)  \mathbb{P}_{k}(A_{i})=\left(  1-\frac{1}{8}a^{4}M\right)
\mathbb{P}_{k}(A).
\end{multline*}
Whence, by (\ref{rep T}),\ we obtain%
\begin{align*}
\int|T_{k-1}\circ T_{k}(h)(x)|\mathbb{P}_{k-2}(dx) &  \leq\int|h(z)|m(dz)\\
&  \leq\left(  1-\frac{1}{8}a^{4}M\right)  \int|h(z)|\mathbb{P}_{k}(dz),
\end{align*}
and Lemma \ref{estimate3} follows by noting that%
\[
M=1-|f_{k-1}(u)|^{2}.
\]
$\square$

\bigskip

\begin{lemma}
\label{estimate2} For any $k\in\mathbb{N}$ and $u\in\mathbb{R}$, we have
\[
||T_{u,k}(1)||_{\infty}^{2}\leq1-a^{2}(1-|\mathbb{E}(\exp(iuX_{k}%
))|^{2})\text{.}%
\]

\end{lemma}

\noindent\textbf{Proof}: We start by noticing that
\begin{align*}
\left\vert T_{u,k}(1)(x)\right\vert ^{2}  &  =\int\! \!\!\int\exp
(iu(g_{k}(y)-g_{k}(y^{\prime})))Q_{k}(x,dy)Q_{k}(x,dy^{\prime})\\
&  = \int\! \!\!\int\cos(u(g_{k}(y)-g_{k}(y^{\prime})))Q_{k}(x,dy)Q_{k}%
(x,dy^{\prime})\\
&  =1-2\int\! \!\!\int\sin^{2} \Big (\frac{u}{2}(g_{k}(y)-g_{k}(y^{\prime})
)\Big )Q_{k}(x,dy)Q_{k}(x,dy^{\prime}).
\end{align*}
By (\ref{main condition}) we have
\begin{multline*}
\int\! \!\!\int\sin^{2} \Big (\frac{u}{2}(g_{k}(y)-g_{k}(y^{\prime}) )
\Big )Q_{k}(x,dy)Q_{k}(x,dy^{\prime})\\
\geq a^{2} \int\! \!\!\int\sin^{2} \Big (\frac{u}{2}(g_{k}(y)-g_{k}(y^{\prime
}) ) \Big ) \mathbb{P}_{k}(dy)\mathbb{P}_{k}(dy^{\prime}).
\end{multline*}
and the result follows. $\square$

\bigskip

We are now in the position to prove Proposition \ref{factor}.

\bigskip

\noindent\textbf{Proof of Proposition \ref{factor}}

\bigskip

Note that, for $k\geq1$%
\[
\mathbb{E}(\exp(iuS_{2k})|\xi_{0})=T_{1}\circ T_{2}\circ\cdots\circ
T_{2k}(1)(\xi_{0}).
\]
So%
\[
\mathbb{E}(\exp(iuS_{2k}))=\int T_{1}\circ T_{2}\circ\cdots\circ
T_{2k}(1)(x)\mathbb{P}_{0}(dx),
\]
and then
\[
|\mathbb{E}(\exp(iuS_{2k}))|\leq||T_{1}\circ T_{2}|| \cdots||T_{2k-1}\circ
T_{2k}||\text{ .}%
\]
Hence, by Lemma \ref{estimate3} we have that, for $k\geq1,$%
\[
|\mathbb{E}(\exp(iuS_{2k}))|\leq\prod_{j=1}^{k}[1-a^{4}(1-|\mathbb{E}%
(\exp(iuX_{2j-1})|^{2})/8]\text{ .}%
\]
Also, by Lemma\ \ref{estimate3} and Lemma \ref{estimate2}
\begin{gather*}
|\mathbb{E}(\exp(iuS_{2k}))|\leq||T_{2}\circ T_{3}|| \cdots||T_{2k-2}\circ
T_{2k-1}||\text{ }||T_{2k}(1)||\\
\leq\prod_{j=1}^{k}[1-a^{4}(1-|\mathbb{E}(\exp(iuX_{2j})|^{2})/8]\text{ ,\ }%
\end{gather*}
and so, by multiplying these two inequalities we get%
\[
|\mathbb{E}(\exp(iuS_{2k}))|^{2}\leq\prod_{j=1}^{2k}[1-a^{4}(1-|\mathbb{E}%
(\exp(iuX_{j})|^{2})/8]\text{ }%
\]
A similar result can be obtain for $|\mathbb{E}(\exp(iuS_{2k+1}))|^{2}$, and
therefore
\[
|{\mathbb{E}}(\exp(iuS_{n}))|^{2}\leq\prod_{j=1}^{n}\Big [1-\frac{a^{4}}%
{8}(1-|f_{j}(u)|^{2})\Big ].
\]
Now, for any real $x$, $1+x\leq\exp(x)$, and the result in Proposition
\ref{factor} follows. $\square$

\subsection{Proof of Theorem \ref{ThLocalAB}}

We consider first the non-lattice case. According to Remark \ref{remd1d2}, it
remains to verify conditions $(D_{1})$ and $(D_{2})$. Suppose $D>0$. Choose
any $\delta$ such that $0<\delta<D$. By Proposition \ref{factor} combined with
Condition \textbf{A}, for any $1\leq|u|<\delta B_{n}$,
\begin{align*}
\Big |{\mathbb{E}}\exp\Big (iu\frac{S_{n}}{B_{n}}\Big )\Big | &  \leq
\exp\Big [-\frac{a^{4}}{16}\sum_{j=1}^{n}(1-|f_{j}(\frac{u}{B_{n}}%
)|^{2})\Big ]\\
&  \leq\exp(-g(|u|)).
\end{align*}
Integrating both sides of this inequality on the intervals $T\leq|u|\leq\delta
B_{n},$ with the restriction $T\geq1$, we obtain
\begin{align*}
\int\nolimits_{T\leq|u|\leq\delta B_{n}}\Big |{\mathbb{E}}\exp\Big (iu\frac
{S_{n}}{B_{n}}\Big )\Big |du &  \leq\int\nolimits_{T\leq|u|\leq\delta B_{n}%
}\exp(-g(|u|))du\\
&  \leq\int\nolimits_{|u|\geq T}\exp(-g(|u|))du.
\end{align*}
Whence, taking first $\limsup_{n}$ and then $T\rightarrow\infty,$ condition
$\left(  D_{1}\right)  $ is verified.

We move now to verify $\left(  D_{2}\right)  .$ Because the interval
$[-D,-\delta]\cup$ $[\delta,D]$ is compact, $\left(  D_{2}\right)  $ is
verified if we can show that for any $u$ fixed in this interval we can find an
open interval $O_{u}$ containing $u$ such that
\begin{equation}
B_{n}\sup_{t\in O_{u}}|{\mathbb{E}}\exp(itS_{n})|\rightarrow0\text{ as
}n\rightarrow\infty.\label{cond Diii}%
\end{equation}
By Proposition \ref{factor}, for any $t$,
\begin{align*}
B_{n}|{\mathbb{E}}\exp(itS_{n})| &  \leq B_{n}\exp\Big [-\frac{a^{4}}{16}%
\sum\nolimits_{j=1}^{n}(1-|f_{j}(t)|^{2})\Big ]\\
&  =\exp\Big [\ln B_{n}-\frac{a^{4}}{16}\sum_{k=1}^{n}(1-|f_{k}(t)|^{2}%
)\Big ].
\end{align*}
Now (\ref{cond Diii}) is satisfied provided that
\[
\ln B_{n}\left(  1-\inf_{t\in O_{u}}\frac{1}{\ln B_{n}}\frac{a^{4}}{16}%
\sum_{k=1}^{n}(1-|f_{k}(t)|^{2})\right)  \rightarrow-\infty.
\]
Since $B_{n}\rightarrow\infty$, we obtain in this case that $(D_{2})$ follows
from Condition \textbf{B}.

The proof of the lattice case is similar.

$\ \square$

\subsection{Proof of Theorem \ref{normal attr} and Remark \ref{remnormal attr}%
}

The proof of Theorem \ref{normal attr} is similar to the proof of the
corresponding result in Merlev\`{e}de et al. (2021) replacing Proposition 8
there by Proposition \ref{factor} here. In case $\mathbb{E}(X_{0}^{2}%
)<\infty,$ the identification of normalizer $B_{n}$ in Remark
\ref{remnormal attr} follows by (\ref{varineq}) and (\ref{limvar}).

When $\mathbb{E}(X_{0}^{2})=\infty,$ we shall invoke Theorem 2.1 in Peligrad
(1990) to decide that $B_{n}$ can be taken $\sqrt{\pi/2}\mathbb{E}|S_{n}|$. To
apply this theorem, note that by (\ref{phi1sm1}) we have $\varphi_{1}<1,$ and
also the fact that $X_{0}$ is in the domain of attraction of the normal law
with $\mathbb{E}(X_{0}^{2})=\infty,$ is equivalent to (1.2) in Peligrad
(1990). Moreover, by Lemma 5.2 and Corollary 5.1 in Peligrad (1990) it follows
that $B_{n}\rightarrow\infty.$

\subsection{Proof of Theorem \ref{Th attraction}}

We shall verify first condition (\ref{conv distribution}). To derive this
convergence in distribution, we use the existing results, inspired by Theorem
3.2 in Samur (1987) and further developed by Kobus (1995) and
Tyran-Kami\'{n}ska, M. (2010). As a matter of fact, we shall apply Corollary
(5.9) in Kobus (1995).

Note that, by relations between the mixing coefficients in Section \ref{mix},
our Markov chain is $\varphi-$mixing with exponential rate of convergence to
$0$ of the $\varphi$-mixing coefficients. Therefore, also the sequence
$(X_{k})_{k\in\mathbb{Z}}$ is $\varphi-$mixing with exponential rate of
convergence to $0$, and also $\rho-$mixing$\ $\ with exponential rate of
convergence to $0.$ Whence the two mixing conditions (i) and (ii) in Corollary
(5.9) in Kobus (1995) are satisfied. Since we assumed that the distribution of
$X_{0}$ satisfies conditions (\ref{tp}) and (\ref{lrtp}), and in addition
$\mathbb{E}(X_{0})=0$ for $1<p<2$ and $X_{0}$ has a symmetric distribution for
$p=1$, by classical results, it follows that $X_{0}$ is in the domain of
attraction of a strictly stable distribution. This gives that condition (5.12)
in Kobus (1995) holds without centering. By the conclusion of Corollary (5.9)
in Kobus (1995) in order to prove
\[
\frac{S_{n}}{B_{n}}\Rightarrow L,
\]
with $L$ strictly stable with exponent $p,$ it is enough to verify the
following condition: for any $x>0$ and any $k\in\mathbb{N}$ we have
\[
\lim_{n\rightarrow\infty}n\mathbb{P}\left(  |X_{0}|\geq xB_{n},|X_{k}|\geq
xB_{n}\right)  =0.
\]
Clearly, by conditions (\ref{tp}) we can write
\[
x^{p}n\mathbb{P}\left(  |X_{0}|\geq xB_{n}\right)  =\frac{n}{B_{n}^{p}}%
\ell(B_{n})\frac{\ell(xB_{n})}{\ell(B_{n})}.
\]
From the properties of slowly varying functions and from (\ref{defBn}) we
deduce that for any $x>0,$%
\[
\lim_{n\rightarrow\infty}n\mathbb{P}\left(  |X_{0}|\geq xB_{n}\right)
=\frac{1}{x^{p}}.
\]
Therefore we can find a constant $C_{x}$ such that%
\begin{align*}
&  n\mathbb{P}\left(  |X_{0}|\geq xB_{n},|X_{j}|\geq xB_{n}\right) \\
&  =n\mathbb{P}\left(  |X_{0}|\geq B_{n}x\right)  \left(  \mathbb{P}%
|X_{j}|\geq xB_{n}\text{ }|\text{ }|X_{0}|\geq xB_{n}\right) \\
&  \leq C_{x}\mathbb{P}\left(  |X_{j}|\geq xB_{n}\text{ }|\text{ }|X_{0}|\geq
xB_{n}\right)  ,
\end{align*}
and the result follows by condition (\ref{DJ}).

Also note that a stable distribution with index $p\neq1$ and symmetric for
$p=1,$ has an integrable characteristic function. This can be seen by the form
of the modulus of the characteristic function for these situations, namely for
a constant $c_{p}$:%
\[
|f_{L}(t)|=\exp(-|c_{p}t|^{p}).
\]
(L\'{e}vy, 1954). Therefore condition (\ref{conv distribution}) holds.

To prove the results, according to Theorems \ref{general} and \ref{ThLattice},
it remains to verify conditions $\mathbf{A}$ and $\mathbf{B}$.

We shall invoke relation (7.9) in Feller (1967), which states that under
(\ref{tp}) and (\ref{lrtp}),\ there are constants $\delta>0$, $n_{0}%
\in\mathbb{N}$ and $c>0$ such that for $0<|u|\leq\delta B_{n}$ and $n>n_{0}$
\[
n\left(  1-|f_{0}(\frac{u}{B_{n}})|^{2}\right)  >c|u|^{p}.
\]
Clearly Condition \textbf{A} is satisfied.

In the stationary setting Condition \textbf{B} reads: For $u\neq0$ there is an
$\varepsilon=\varepsilon(u)$, $c(u)$ and a $n_{0}=n_{0}(u)$ such that for all
$t$ with $|t-u|\leq\varepsilon$ and $n>n_{0},$
\[
\frac{na^{4}}{2^{4}(\ln B_{n})}(1-|f_{0}(t)|^{2})\geq c(u)>1.
\]
On another hand, by the definition of $B_{n}$, and the properties of slowly
varying functions, it is well known that $B_{n}<n^{p^{\prime}}$ for
$p^{\prime}>1/p$ and $n$ large enough. So $\ln B_{n}<p^{\prime}\ln n$ for
$p^{\prime}>1/p$ and $n$ sufficienly large$.$

Therefore, if $X$ does not have a lattice distribution, then for every
$u\neq0,$ $|f_{0}(u)|<1,$ and by the continuity of $f_{0}$, for any $u\neq0$
we can certainly find an $\varepsilon>0$ and $d<1$, such that for all $t$ with
$|t-u|\leq\varepsilon$ we have $|f_{0}(t)|<d.$ Therefore, in this case, for
such $t$,
\[
\frac{na^{4}}{(\ln B_{n})}(1-|f_{0}(t)|^{2})\geq\frac{na^{4}}{p^{\prime}\ln
n}(1-d^{2})\rightarrow\infty
\]
and Condition \textbf{B} is satisfied.

In the lattice case we use a similar argument on the interval $0<u<\frac{\pi
}{h}$. $\ \square$

\section{Acknowledgement}

The authors are grateful to the referee for carefully reading the paper and
for insightful corrections and suggestions that significantly improved the
presentation of the paper. This paper was partially supported by the NSF grant DMS-2054598.

\end{document}